\documentclass[12pt,reqno]{article}

\usepackage[usenames]{color}

\usepackage{graphicx}
\usepackage{amscd}
\usepackage{array}

\definecolor{webgreen}{rgb}{0,.5,0}
\definecolor{webbrown}{rgb}{.6,0,0}

\usepackage{fullpage}
\usepackage{float}

\usepackage{graphics,amsmath,amssymb}
\usepackage{amsthm}
\usepackage{amsfonts}
\usepackage{latexsym}
\usepackage{epsf}
\usepackage{mathrsfs}

\usepackage[colorlinks=true,
linkcolor=webgreen,
filecolor=webbrown,
citecolor=webgreen]{hyperref}

\newcommand{\seqnum}[1]{\href{http://oeis.org/#1}{\underline{#1}}}

\def \Leg{\overwithdelims()}

\begin{document}

\theoremstyle{plain}
\newtheorem{theorem}{Theorem}
\newtheorem{corollary}[theorem]{Corollary}
\newtheorem{lemma}[theorem]{Lemma}
\newtheorem{proposition}[theorem]{Proposition}

\theoremstyle{definition}
\newtheorem{definition}[theorem]{Definition}
\newtheorem{example}[theorem]{Example}
\newtheorem{conjecture}[theorem]{Conjecture}

\theoremstyle{remark}
\newtheorem*{remark}{Remark}

\begin{center}
\vskip 1cm{\LARGE\bf  Periodic Weighted Sums of Binomial Coefficients}

\vskip 1cm
\large

Greg Dresden\\ 
Washington \& Lee University \\
Lexington, VA, 24450 \\ USA\\
\href{mailto:dresdeng@wlu.edu}{\tt dresdeng@wlu.edu} \\

\ \\

Yike Li\\
Beijing City International School\\
Shuangjing, Chaoyang, Beijing, 100022\\ China\\
\href{redcedariver@hotmail.com}{\tt redcedariver@hotmail.com}\\ 
\end{center}
\vskip0.2in

\begin{abstract} 
Using elementary methods, we establish old and new relations between binomial coefficients, Fibonacci numbers, Catalan numbers, Pell numbers, and more.
\end{abstract}

\section{Introduction}
We present a general theorem on linear recurrences and we use it to produce  results both old and new about  weighted sums of  binomial coefficients. One highlight is this first equation, which connects the Fibonacci numbers, the binomial coefficients, and the Legendre symbol:
\begin{equation}\label{e.first}
F_{2n} = 
\sum_{k=0}^n {k \Leg 5}{2n \choose n+k}
 = 
    \frac{1}{5^n} \sum_{k=0}^{2n}
    (-1)^{k+1}{k \Leg 5} {4n \choose 2n+k}.
\end{equation}
The first equality can also be derived  from a more complicated expression of Andrews \cite{Andrews} from 1969, but the second equality is new.

Also of interest are the non-negative integer solutions $(X_n, Y_n)$ 
to the Pell equation $X^2 - 3Y^2 = 1$. We show that 
\begin{equation}\label{e.Xn}
2X_n +2^n = 3\sum_{j=-n}^n (-1)^j {2n \choose n+6j}   
\end{equation}
and also that 
\begin{equation}\label{e.Yn}
Y_n = {2n \choose n+1} 
    - {2n \choose n+5} - {2n \choose n+7}
    + {2n \choose n+11} + {2n \choose n+13}
    - \cdots
\end{equation}              
and while the formula for $X_n$ follows from Merca \cite{Merca}, the formula for $Y_n$ is new. 

To provide some  background, 
we note that we began this research project by looking at the Binet formula for $F_{2n}$, which is
\begin{equation}\label{e.Binet2n}
F_{2n} = \frac{1}{\sqrt{5}}\left(
\left(
\frac{1+\sqrt{5}}{2}
\right)^{2n}
-
\left(\frac{1-\sqrt{5}}{2}\right)^{2n}
\right).
\end{equation}
Since 
$ (1+\sqrt{5})/2 = 2\cos \pi/5$ and 
$(1-\sqrt{5})/2 = 2\cos 3\pi/5$, then the above equation becomes
\begin{equation}\label{e.F2ncos}
F_{2n} = \frac{1}{\sqrt{5}}(
( 2 \cos \pi/5 )^{2n}
-
( 2 \cos 3\pi/5 )^{2n}
),
\end{equation}
and we then realized that we could 
apply a cosine formula to expand the above equation into a more interesting form, which eventually gave us the first part of equation (\ref{e.first}). 
After applying this same technique to other sequences such as $(X_n)_{n\geq 0}$ and $(Y_n)_{n\geq 0}$ mentioned earlier, we soon recognized that the Binet formula for these sequences was not necessary; in each case we only needed the linear recurrence formula itself. 

In what follows, we will establish a general theorem which covers a large class of linear recurrence sequences. To be precise, 
 so long as the coefficients of the linear recurrence satisfy an easily-verified property, then  our theorem 
 will give us an identity involving a weighted sum of binomial coefficients as seen in equations 
(\ref{e.first}), (\ref{e.Xn}), and (\ref{e.Yn}) in which the weights,
typically $1$ and $-1$, appear in a pattern (hence the title of our paper). Sometimes this periodic pattern of weights matches up with a Legendre or Kronecker symbol as seen in equation (\ref{e.first}), but other times 
as seen in equation (\ref{e.Yn}) the best way to express the pattern is to write it out directly.  Furthermore, the values of these periodic weights   depend  {\em only} on the initial values of the linear recurrence sequence.
We  then  use this theorem to derive new identities for the Fibonacci numbers, the Pell numbers, the Catalan numbers, and more. Many of our identities now appear in the Online Encyclopedia of Integer Sequences (OEIS) \cite{oeis}.

\section{Main result}

To motivate some definitions, let us return back to our numbers $F_{2n}$, which satisfy the  equation
\begin{equation}\label{e.F2n3neg1}
F_{2n} = 3 F_{2(n-1)} -1 F_{2(n-2)} \qquad \mbox{for $n \geq 2$}.
\end{equation}
We say that $(F_{2n})_{n\geq 0}$ is a 
{\em linear recurrence sequence}, and 
since each term in equation (\ref{e.F2n3neg1}) is defined by its two previous terms, we say that the sequence has {\em order} 2. 
 The set of coefficients on the right of 
equation (\ref{e.F2n3neg1}), in this case $\{3, -1\}$, is called the {\em signature} of this linear recurrence. 
 From equation (\ref{e.F2n3neg1}) we now define the 
{\em characteristic polynomial} for the 
sequence  $(F_{2n})_{n \geq 0}$
to be 
\[
x^2 - (3x - 1).
\]
Of course, all these definitions generalize nicely to linear recurrences of arbitrary order. 

Now, one of the roots of our  characteristic polynomial is
$((1 + \sqrt{5})/2)^2 = (2 \cos \pi/5)^2$. Another way of expressing this is to say that $x^2-(3x-1)$ is  the {\em minimal polynomial} for $(2 \cos \pi/5)^2$. 

With this in mind, we are ready to present our main theorem.

\begin{theorem}\label{c.1} 
Suppose we have a linear recurrence sequence $(A_n)_{n \geq 0}$  
such that its characteristic polynomial
is the minimal polynomial for  $(2 \cos \pi/q)^2$ for some  integer $q\geq 3$. Then, we have
\begin{equation}\label{e.An2}
A_n = A_0 {2n \choose n} + \sum_{k=1}^n \omega_k {2n \choose n+k},
\end{equation}
where we define
$\omega_k$ in terms of $A_0, A_1, A_2, \dots$ as follows: $\omega_0 = 2A_0$, and 
\begin{equation}\label{e.wk2}
\omega_k = 
\sum_{j=0}^k (-1)^j{2k-j \choose j} \frac{2k}{2k - j}A_{k-j} \qquad \mbox{for $k \geq 1$}.
\end{equation}
Furthermore, 
 we have that     $\omega_{k'}  =  \omega_{k}$ 
for $k' \equiv \pm k$  {\rm (}mod $q${\rm )},
and if $q$ is even then $\omega_{k'} = -\omega_{k}$ 
for $k' \equiv q/2 - k$  {\rm (}mod $q${\rm )}.
\end{theorem}
For  the proof, see Section \ref{s.proof} at the end of this article. 

We will find it useful to write out the first few values of $\omega_k$. We use  equation (\ref{e.wk2}), along with our definition of $\omega_0$ as $2A_0$, to get the values in Table \ref{e.omegak012}.

\begin{table}[h!]
\centering
\begin{tabular}{l} 
$\omega_0 = 2A_0$\\[1.0ex]
$\omega_1 = A_1 - 2A_0$ \\[1.0ex]
$\omega_2 = A_2 - 4A_1 + 2A_0$ \\[1.0ex]
$\omega_3 = A_3 - 6A_2 + 9A_1  - 2A_0$ \\
\end{tabular}
\caption{Initial values of $\omega_k$ from Theorem \ref{c.1}.}\label{e.omegak012}
\end{table}

The coefficients of $A_n$ in Table \ref{e.omegak012} appear in reverse order as sequence \seqnum{A127677}
in the OEIS.

Now, Theorem \ref{c.1} applies to 
a large class of linear recurrences. 
One such example that we will see in Section \ref{s.Catalan} is the third-order sequence
\[
A_n = 5A_{n-1} - 6A_{n-2} 
+ 1A_{n-3}
\]
with signature $\{5, -6, 1\}$ and characteristic polynomial 
$x^3 - (5x^2 -6x+1)$, because 
that polynomial 
is the minimal polynomial for $(2 \cos \pi/q)^2$ with $q=7$.
For convenience, we list in 
Table \ref{tab.1} all such cases for linear recurrences with order no more than three, indexed by the corresponding values of $q$.

\renewcommand{\arraystretch}{1.3}
\begin{table}[H]
\centering
\begin{tabular}{r|lcr|l} 
 $q$ & Signature & \ \ \ \ \ \ \ \ \ \ \ & $q$ & Signature \\ 
 \hline
 5 &$ \{3,-1\} $  &  & 10 & $\{5,-5\} $\\ 
 7 & $\{5,-6,1\} $ &   & 12 & $\{4,-1\} $\\ 
 8 & $\{4,-2\} $  &    & 14 & $\{7,-14,7\} $\\ 
 9 & $\{6,-9,1\} $ &   & 18 & $\{6,-9,3\}$
\end{tabular}
 \caption{All signatures of length 2 or 3 for minimal polynomials for $(2 \cos \pi/q)^2$.}
\label{tab.1}
\end{table}

Thus, for any entry in the OEIS with one of these signatures, we can write the terms in that sequence as a weighted sum of binomial coefficients as given in equation (\ref{e.An2}) from Theorem \ref{c.1}. We can do the same for signatures of length 4 or more, but those are not as interesting. Let us now look at some examples.

\section{Applications}\label{S2}

\subsection{The Fibonacci numbers}

As mentioned above, the  sequence $(F_{2n})_{n\geq 0}$ 
has signature $\{3, -1\}$ which is the first entry in our Table \ref{tab.1}.
This gives us the following formula for $F_{2n}$ as a weighted sum of binomial coefficients.

\begin{theorem}\label{t.31}
For $F_{n}$ the Fibonacci numbers, we have 
\begin{equation}
F_{2n}  = 
\sum_{k=0}^n {2n \choose n+k}
{k \Leg 5},\label{e.F2n1}
\end{equation}
where   $\displaystyle {k \Leg 5}$ represents the Legendre symbol.
\end{theorem}

\begin{remark}
As pointed out by our helpful referee, the appearance of the Legendre symbol in Theorem 
\ref{t.31} should not come as a complete surprise. If we were to 
apply the fairly well-known cosine equation (\ref{e.cosine2}) directly to our Binet formula (\ref{e.F2ncos}) for $F_{2n}$, we 
would get 
\begin{equation}\label{e.F2ncosNEW}
F_{2n} = \frac{1}{\sqrt{5}}\sum_{k=1}^n {2n \choose n+k} 
\left( 2 \cos 2k\pi/5 
-
 2 \cos 6k\pi/5 \right).
\end{equation}
For $k$ a multiple of 5, then 
the expression $(2 \cos 2k\pi/5 
-
 2 \cos 6k\pi/5)$ equals $0$, but otherwise it has the same values as the 
quadratic Gauss sum 
\[
g(k;5) = \sum_{n=0}^4 e^{2 \pi i k n^2/5} .
\]
From a  theorem by Berndt, Evans, and Williams \cite[Theorem 1.5.2]{BEW}, we learn that this satisfies
\[
g(k;5) = {k \Leg 5} \sqrt{5} \qquad \mbox{for $k$ not a multiple of 5,}
\]
and so using this in equation (\ref{e.F2ncosNEW}) would give us a direct proof of our Theorem \ref{t.31}.
\end{remark}

\begin{proof}[Proof of Theorem \ref{t.31}]
If we set $A_n = F_{2n}$ then $A_n$ satisfies the recurrence $A_n = 3A_{n-1} - A_{n-2}$. Thanks to  Table \ref{tab.1}, 
we can apply  Theorem \ref{c.1}
with  
$q=5$ and with initial values $A_0 = 0$ (because $F_{2n} =0$ at $n=0$) and $A_1 = 1$ (because $F_{2n} = 1$ at $n=1$). 
From Table \ref{e.omegak012} we have the following values for 
$\omega_k$: 
		\begin{align*}
			\omega_0 &= 2A_0 = 0,\\
			  \omega_1 &= A_1 - 2A_0 = 1 \\ 
            \omega_2 &= A_2 - 4A_1 + 2A_0 = -1 
		\end{align*}
Theorem \ref{c.1} also tells us  that $\omega_k = \omega_{5-k} = \omega_{5+k}$. Using this, along with the three values already given above, we conclude that 
 $\omega_k = 1$ for $k \equiv 1,4$ (mod 5), and  
 $\omega_k = -1$ for $k \equiv 2,3$ (mod 5), 
 and  
 $\omega_k = 0$ for $k \equiv 0$ (mod 5). This means that 
$\omega_k$ has the same values as the Legendre symbol ${k \Leg 5}$, 
and since we also have $A_0 = 0$ then equation (\ref{e.An2})
of Theorem \ref{c.1} gives us our desired equation 
(\ref{e.F2n1}) for $A_n = F_{2n}$. 
\end{proof}

As we mentioned above, our equation (\ref{e.F2n1}) 
for the Fibonacci numbers
is not entirely new.  Andrews \cite{Andrews} used complex numbers to show that 
\[
F_n = \sum_{\alpha = -\infty}^\infty (-1)^{\alpha} {n \choose \lfloor (n-1-5\alpha)/2\rfloor},
\]
where $\lfloor\ \ \rfloor$ represents the greatest integer function. From this, we can obtain our equation (\ref{e.F2n1}) with a bit of work. However, our method is both more general and more direct.

\subsection{Binomial transforms of the Pell and Pell-Lucas numbers}\label{s.bin}

The sequence $(F_{2n})_{n \geq 0}$ from the previous section has the nice property that it is the {\em binomial transform} of the ``regular" Fibonacci sequence. In other words,
\[
F_{2n} = \sum_{i=0}^n {n \choose i} F_i.
\] 
The following theorem will allow us to consider two other sequences of numbers
that are also binomial transforms.

\begin{theorem}\label{t.4neg2}
For $A_n = 4 A_{n-1} - 2A_{n-2}$
 with initial values $A_0$ and $A_1$, then 
\begin{equation}
    A_n = A_0 {2n \choose n} +  
            \sum_{k = 1}^n \omega_k  {2n \choose n+k},
\end{equation} 
with $\omega_k$ repeating modulo 8 as given 
below in Table \ref{tab.w8}, with $\omega_0 =2A_0$ and 
$\omega_1 =A_1 - 2A_0$. 
\end{theorem}

\begin{table}[H]
\centering
\begin{tabular}{c|ccccccccc} 
  $k$ (mod 8) & 0&1 & 2 &\ \ 3&\ \ 4&\ \ 5& 6&7& 8\\ 
  \hline
  $\omega_k$ & $\omega_0$ & $\omega_1$ & 0 & $-\omega_1$ & $-\omega_0$ & $-\omega_1$ & 0 & $\omega_1$ & $\omega_0$
\end{tabular}
 \caption{Values of $\omega_k$ for signature $\{4,-2\}$}
\label{tab.w8}
\end{table}

\begin{proof}
Since the signature $\{4, -2\}$ appears in Table \ref{tab.1},    we can apply Theorem \ref{c.1} with  $q=8$. The values for $\omega_0$ and 
$\omega_1$  come to us from  Table \ref{e.omegak012}, and the remaining  values of $\omega_k$ follow from the conclusion of Theorem    \ref{c.1} 
which tells us that 
    $\omega_k = \omega_{8 \pm k}$ and $\omega_k = -\omega_{4-k}$. 
\end{proof}

With this theorem in hand, we can now produce two  results, one new and one old, about the binomial transforms of the Pell and the Pell-Lucas numbers. 

\begin{corollary}
    For $P_n$ the Pell numbers 
    $0, 1, 2, 5, 12, 29, \dots$ from
    {\rm \seqnum{A000129}}, 
    then
\begin{equation}\label{e.12}
    \sum_{i=0}^n {n \choose i} P_i
  = \sum_{k=0}^n \left(\frac{k}{8}\right)  {2n \choose n+k}
\end{equation}
where   $\displaystyle {k \Leg 8}$ represents the Kronecker symbol.
\end{corollary}
\begin{proof}
If we set $A_n$ equal to the binomial transform of $P_n$ as seen on the left of equation (\ref{e.12}), then the sequence $(A_n)_{n \geq 0}$ begins with 
$0, 1, 4, 14, 48, 164, \dots$ and is
given by \seqnum{A007070} where we also learn that it has signature $\{4, -2\}$. Hence, we can apply 
Theorem \ref{t.4neg2} with $\omega_0 = 2A_0 = 0$ and $\omega_1 = A_1 - 2A_0 
= 1$. Furthermore, Table \ref{tab.w8} from
Theorem \ref{t.4neg2} becomes

\begin{center}
\begin{tabular}{c|ccccccccc} 
  $k$ (mod 8) & 0 & 1 & 2 & \ \ 3& 4&\ \  5& 6&7& 8\\ 
  \hline
  $\omega_k$  & 0 & 1 & 0 & $-1$ & $0$ & $-1$ & 0 & 1 & 0
\end{tabular}
\end{center}

and from this we see that $\omega_k = 
\displaystyle {k \Leg 8}$, as desired.
\end{proof}

\begin{corollary}\label{c.5}
    For $Q_n$ 
    the  Pell-Lucas numbers
    $2, 2, 6, 14, 34, 81, \dots$ from
    {\rm \seqnum{A002203}}, 
    then
\begin{equation}\label{e.13}
    \sum_{i=0}^n {n \choose i} Q_i
= 2{2n \choose n}  + 4\sum_{j \geq 1} (-1)^j {2n \choose n+4j}.
\end{equation}
\end{corollary}
\begin{proof}
The binomial transform of $Q_n$ is 
$2,4,12,40,136, \dots$, which is the sequence 
\seqnum{A056236} with signature $\{4, -2\}$. Hence, if we let 
 $A_n$ equal  the binomial transform of $Q_n$ as seen on the left of equation (\ref{e.13}), then
 we can again apply 
Theorem \ref{t.4neg2}, 
this time with $\omega_0 = 2A_0 = 4$ and $\omega_1 = A_1 - 2A_0 
= 0$. Thus, Table \ref{tab.w8} from
Theorem \ref{t.4neg2}, starting at $k=1$, becomes

\begin{center}
\begin{tabular}{c|cccccccc} 
  $k$ (mod 8) & 1 & 2 &  3& \ \ 4&  5& 6&7& 8\\ 
  \hline
  $\omega_k$  & 0 & 0 & $0$ & $-4$ & $0$ & 0 & 0 & $4$
\end{tabular}
\end{center}
and from this we obtain our desired formula.
\end{proof}

We note that Corollary \ref{c.5} is not a new result.
From Merca
\cite[Corollary 8]{Merca} we have the formula
\[
(2+\sqrt{2})^n + (2-\sqrt{2})^n 
= 2{2n \choose n}  + 4\sum_{j \geq 1} (-1)^j {2n \choose n+4j},
\]
and as seen in \seqnum{A056236} the binomial 
transform of $Q_n$ is indeed equal to 
$(2+\sqrt{2})^n + (2-\sqrt{2})^n$.

As a special case of Corollary \ref{c.5}, we have 
the following relationship between a weighted sum of binomial coefficients (on the left) and a periodic weighted sum of binomial coefficients (on the right). 

\begin{corollary}\label{c.6}
    For $Q_n$ 
    the  Pell-Lucas numbers
    $2, 2, 6, 14, 34, 81, \dots$ as
given by 
    {\rm \seqnum{A002203}}, 
    then
\begin{equation}\label{e.14}
    \sum_{i=0}^{4n} {4n \choose i} Q_i
= 2(-1)^n\sum_{j=0}^{2n} (-1)^j {8n \choose 4j}.
\end{equation}
\end{corollary}
\begin{proof}
    We can split up the right-hand side of equation (\ref{e.14}) as follows:
\[
 2 (-1)^n\sum_{j=0}^{n-1} (-1)^j {8n \choose 4j}
+ 
2{8n \choose 4n} 
+ 
2 (-1)^n\sum_{j=n+1}^{2n} (-1)^j {8n \choose 4j}.
\]
In the first sum we replace $j$ with $n-j$, and in the second sum we replace $j$ with $n+j$, giving us
\[
 2 (-1)^n\sum_{j=1}^{n} (-1)^{n-j} {8n \choose 4n-4j}
+ 
2{8n \choose 4n} 
+ 
2 (-1)^n\sum_{j=1}^{n} (-1)^{n+j} {8n \choose 4n+4j}.
\]
We now use the symmetry of the binomial coefficients to combine the two sums to obtain
\[
2{8n \choose 4n} 
+ 
4 \sum_{j\geq 1} (-1)^{j} {8n \choose 4n+4j},
\]
and by equation (\ref{e.13}) in Corollary \ref{c.5} this is equal to 
$  \sum_{i=0}^{4n} {4n \choose i} Q_i$, as desired. 
\end{proof}

\subsection{Solutions to Pell's Equation}

The non-negative integer solutions $(X_n, Y_n)$ to  the Pell equation $X^2 - 3Y^2 = 1$
are given by 
\seqnum{A001075} for $X_n$, and \seqnum{A001353} for $Y_n$, where 
the first sequence begins $1, 2, 7, 26, 97, \dots$,
and
the second is $0, 1, 4, 15, 56, \dots$.
These are well-known linear recurrence sequences, and 
both of them have signature $\{4,-1\}$ which corresponds
to $q=12$ in Table \ref{tab.1}. With this in mind, 
we present the following theorem and then 
we show how it gives us the equations 
(\ref{e.Xn}) and (\ref{e.Yn}) for $X_n$ and $Y_n$ as mentioned in  the 
introduction to this paper. 

\begin{theorem}\label{t.4neg1}
For $A_n = 4 A_{n-1} - A_{n-2}$
 with initial values $A_0$ and $A_1$, then 
\begin{equation}\label{e.4neg1}
    A_n = A_0 {2n \choose n} +  
            \sum_{k=1}^n \omega_k  {2n \choose n+k},
\end{equation} 
with $\omega_k$ repeating modulo 12 as given below in Table \ref{tab.w12}, with 
$\omega_1 =A_1 - 2A_0$
and 
$\omega_2 =A_0$.

\begin{table}[H]
\centering
\begin{tabular}{c|ccccccccccccc} 
  $k$ (mod 12) & 0 & 1 & 2 & 3 & 4 & 5 & 6 & 7 & 
  8 & 9 & 10 & 11 & 12 \\ 
  \hline
  $\omega_k$  & $2\omega_2$  & 
  $\omega_1$ & $\omega_2$ & 0 & $-\omega_2$ 
  & $-\omega_1$ & $-2\omega_2$ & $-\omega_1$ & $-\omega_2$ & 
  0 & 
  $\omega_2$ & $\omega_1$  & $2\omega_2$
\end{tabular}
 \caption{Values of $\omega_k$ for signature $\{4,-1\}$}
\label{tab.w12}
\end{table}
\end{theorem}

\begin{proof}
Thanks to Table \ref{tab.1}, we can apply 
Theorem \ref{c.1} with $q=12$. 
The values  $\omega_0=2A_0$ and 
$\omega_1 = A_1 - 2A_0$ are given to us in Table 
\ref{e.omegak012}. As for $\omega_2$, we 
see in Table \ref{e.omegak012} that
    $\omega_2 = A_2 - 4A_1 + 2A_0$, but since $A_2 = 4A_1 - A_0$ then this becomes just $\omega_2 = A_0$. Hence, we can write $\omega_0 = 2 \omega_2$. The remaining values of $\omega_k$ follow from the conclusions of Theorem    \ref{c.1} that $\omega_k = \omega_{12 \pm k} = -\omega_{6-k}$, allowing us to fill in the rest of Table \ref{tab.w12}.
\end{proof}

Thanks to  Theorem \ref{t.4neg1}, we have the following corollary. 

\begin{corollary}\label{c.8}
For $X_n$ and $Y_n$ the non-negative solutions to $X^2 - 3Y^2 = 1$, then $X_n$ satisfies equation (\ref{e.Xn}) and $Y_n$ satisfies equation (\ref{e.Yn}).     
\end{corollary}

As we mentioned in the introduction, equation (\ref{e.Xn}) for $X_n$ comes from Merca \cite{Merca}; to be precise, it follows from his Theorem 3 with $n=6$. However, we believe that our equation (\ref{e.Yn}) for $Y_n$ is new. 

\begin{proof}[Proof of Corollary \ref{c.8}]
We begin with $X_n$, which has $X_0 = 1$ and $X_1 =2$. From Theorem \ref{t.4neg1}, we have  $\omega_1 = 0$ and  $\omega_2 = 1$ and so Table \ref{tab.w12} gives us the following values:
\begin{center}
\begin{tabular}{c|ccccccccccccc } 
  $k$ (mod 12) & 0&1 & 2 &3&\ \ 4&5&\ \ 6&7&\ \ 8&9&10&11&12\\ 
  \hline
  $\omega_k$ & 2 & 0 & 1 & 0 & $-1$ & 0 & $-2$ & 0 & $-1$ & 0 & 1 & 0 & 2
\end{tabular}
\end{center}
To reveal a hidden pattern, we write $\omega_k$ as the sum of two periodic sequences, as shown here. We have dropped some of the $0$'s for legibility.
\begin{center}
\begin{tabular}{c|ccccccccccccc } 
  $k$ (mod 12) & \ \ 0&1 & 2 &3&\ \ 4&5&\ \ 6&7&\ \ 8&9&10&11& \ \ 12\\ 
  \hline
    $\omega_k$ & \ \ 2 & 0 & 1 & 0 & $-1$ & 0 & $-2$ & 0 & $-1$ & 0 & 1 & 0 & \ \ $2$\\
  \hline
first sequence & $-1$ &  & 1 &  & $-1$ &  & \ \ $1$ &  & $-1$ &  & 1 & & $-1$\\
second sequence & \ \ $3$ &  &  &  &  &  & $-3$ &  & &  &  & & \ \ 3
\end{tabular}
\end{center}
Thus, equation (\ref{e.4neg1})
tells us that 
\begin{equation}\label{e.15}
 X_n =  {2n \choose n} - 
            \sum_{i=1}^n (-1)^i  {2n \choose n+2j} 
        +  3\sum_{j=1}^n (-1)^j  {2n \choose n+6j}. 
\end{equation}
From Lewis \cite[p.~200]{Lewis}, after adjusting for his notation  
we have the following formula:
\[
2^{2n} \sin^{2n} z = 
{2n \choose n} +  \sum_{h=1}^n (-1)^h  {2n \choose n+h} 2\cos 2hz, 
\]
and if we let $z = \pi/4$, divide by 2, and then renumber the sum, we obtain
\[
2^{n-1}  = 
\frac{1}{2} {2n \choose n} +  \sum_{i=1}^n (-1)^i  {2n \choose n+2i}.
\]
(We could also obtain this from Merca \cite[Corollary 8]{Merca}). When we substitute this into equation (\ref{e.15}) and multiply by two,
we get 
\[
2X_n = 3 {2n \choose n} 
    - 2^{n}
     +  3\sum_{j=1}^n (-1)^j \cdot 2{2n \choose n+6j}.
\]
Since $2{2n \choose n+6j} = {2n \choose n-6j} + {2n \choose n+6j}$, we can re-write that sum on the right to include the terms for $j<0$,
and as for the $j=0$ term we simply bring 
$3 {2n \choose n} $ into the sum. After simplifying,  this gives us equation (\ref{e.Xn}) at the beginning of this paper.

As for $Y_n$, we have  $Y_0 = 0$ and $Y_1 = 1$ and so Theorem \ref{t.4neg1} tells us that in this case we have $\omega_1 = 1$ and $\omega_2 = 0$. 
Thanks to Table \ref{tab.w12}, we quickly assemble the following chart.

\begin{center}
\begin{tabular}{c|ccccccccccccc } 
  $k$ (mod 12) & 0&1 & 2 &3&4&\ \ 5&6&\ \ 7&8&9&10&11 & 12\\ 
  \hline
  $\omega_k$ & 0 & 1 & 0 & 0 & 0 & $-1$ & 0 & $-1$  &   0 & 0 & 0 & 1 & 0
\end{tabular}
\end{center}
From this we immediately get the expression for $Y_n$ in equation (\ref{e.Yn}) from the introduction.
\end{proof}

\subsection{Catalan numbers and a linear recurrence of order three}\label{s.Catalan}

For a particular class of 
linear recurrences, we have the
following theorem. 

\begin{theorem}\label{t.561}
For $A_n$ a sequence satisfying $A_n = 5 A_{n-1} - 6 A_{n-2}+ A_{n-3}$
and with initial values $A_0, A_1$, and $A_2$, then 
\begin{equation}\label{e.561}
    A_n = A_0 {2n \choose n} +  
            \sum_{k =1}^n \omega_k  {2n \choose n+k},
\end{equation} 
with
\begin{equation}
\omega_i = \begin{cases}
2A_0, & \text{for $k \equiv 0$ {\rm (}mod $7${\rm )};}\\[1ex]
A_1 - 2A_0, & \text{for $k \equiv 1,6$ {\rm (}mod $7${\rm )};}\\[1ex]
A_2 - 4A_1 + 2A_0, & \text{for $k \equiv 2,5$ {\rm (}mod $7${\rm )};}\\[1ex]
3A_1 - A_0 - A_2, & \text{for $k \equiv 3,4$ {\rm (}mod $7${\rm )}.}\\[1ex]
\end{cases}\label{e.562}
\end{equation}
\end{theorem}

\begin{proof}
Since the signature $\{5, -6, 1\}$ appears in Table \ref{tab.1},    we can apply Theorem \ref{c.1}, this time  with $q=7$. 
The values for $\omega_0$, $\omega_1$, and $\omega_2$ are given to us in Table \ref{e.omegak012}, and from that same 
equation we also have
\[
\omega_3 = A_3 - 6A_2 + 9A_1  - 2A_0,
\]
but since $A_3 = 5A_2 -6A_1 + A_0$ then this becomes 
\[
\omega_3 = - A_2 + 3A_1  - A_0,
\]
as desired. 
Finally, Theorem    \ref{c.1} also tells us that     $\omega_k = \omega_{7 \pm k}$, thus 
concluding our proof. 
\end{proof}

With this in hand, we can present two new results related to Catalan paths.

\begin{corollary}\label{c.Bn}
    For $B_n$ the  numbers 
    $	1, 1, 2, 5, 14, 42, \dots$ from 
    {\rm \seqnum{A080937}}, 
    then
\begin{equation}\label{e.Bn}
    B_n
= \frac{1}{n+1} {2n \choose n}
 + \sum_{j \geq 1}   4{2n \choose n+7j} - \sum_{j \geq 1}  {2n+2 \choose n+7j+1}
\end{equation}
where   the first term on the right is also known as the $n$th Catalan number {\rm (\seqnum{A000108})}. 
\end{corollary}
\begin{proof}
From \seqnum{A080937} we learn that the numbers $B_n$ have the recurrence
$B_n = 5B_{n-1} - 6B_{n-2} + B_{n-3}$, and this signature $\{5,-6,1\}$ appears in Table \ref{tab.1} with $q=7$. Thus, we can once again apply Theorem \ref{c.1} (this time with $q=7$) and  
from Table \ref{e.omegak012} we 
learn that $\omega_0$, $\omega_1$, $\omega_2$, and $\omega_3$ are $2, -1, 0, 0$ respectively. 
 The remaining values of $\omega_k$ follow from the     fact that $\omega_k = \omega_{7 \pm k}$.
 Here are those values for $\omega_k$, starting at $k=0$.
\begin{center}
\begin{tabular}{c|cccccccc} 
  $k$ (mod 7) & \ \ 0 & 1 & 2 &  3& 4&  5& \ \ 6&7\\ 
  \hline
  $\omega_k$  & 2 & $-1$ & 0 & $0$ & 0 & $0$ & $-1$ & 2 
\end{tabular}
\end{center}
 
 Thus, equation (\ref{e.An2}) from Theorem 
 \ref{c.1} tells us that 
 \begin{equation*}
    B_n
=  {2n \choose n} 
        - {2n \choose n+1}  
 + \sum_{j \geq 1}  
 \left(2 {2n \choose n+7j} - {2n \choose n+7j-1} - {2n \choose n+7j+1} \right).
\end{equation*}
A simple calculation verifies that 
${2n \choose n} - {2n \choose n+1}$ is equal to the Catalan number 
$\frac{1}{n+1}{2n \choose n}$, and 
as for the expression inside the sum, 
we can easily verify that 
\begin{equation*}
2 {2n \choose n+7j} - {2n \choose n+7j-1} - {2n \choose n+7j+1} 
= 4 {2n \choose n+7j} - {2n+2 \choose n+7j+1}.
\end{equation*} 
This gives us equation (\ref{e.Bn}) for $B_n$, as desired. 
\end{proof}

From Corollary \ref{c.Bn} we can obtain an unexpected equation for a type of Catalan paths. We recall \cite[p.~152]{KoshyCat} that the $n$th Catalan number (call it $C_n$) counts the total number of non-negative paths of length $2n$ that start and end at 0, with each step $\pm 1$. Our numbers $B_n$ from Corollary \ref{c.Bn}, as described in \seqnum{A080937}, count the total number of such paths that have maximum height 5 or less. Hence, if we re-write equation (\ref{e.Bn}), recalling again that 
$C_n = \frac{1}{n+1}{2n \choose n}$, we obtain the following.
\begin{corollary}\label{c.cat}
    The number of Catalan paths of length $2n$ with maximum height at least 6 is 
\begin{equation}\label{e.Bn3}
\sum_{j \geq 1} \left( {2n+2 \choose n+7j+1}  -    4{2n \choose n+7j}\right).
\end{equation}
\end{corollary}
This formula gives us the (new) sequence \seqnum{A359311}, which begins 
$0, 0, 0, 0, 0, 0, 1, 12, 89, \dots$. 
We can only imagine that there is an easy combinatorial proof of our 
Corollary \ref{c.cat}, although we are not clever enough to find it. 

\subsection{More binomial transforms}

As we saw in Section \ref{s.bin}, 
the binomial transforms of certain sequences (in that case, the Pell and Pell-Lucas sequences) can be expressed as
weighted sums of binomial coefficients. 
We can do the same for the binomial transform of $F_{2n}$ and $L_{2n}$. These transforms are given in 
\seqnum{A093131} and \seqnum{A020876} respectively, and both have signature
$\{5, -5\}$. This leads us to the following theorem and its three corollaries.

\begin{theorem}\label{t.5neg5}
For $A_n = 5 A_{n-1} - 5A_{n-2}$
 with initial values $A_0$ and $A_1$, then 
\begin{equation}
    A_n = A_0 {2n \choose n} +  
            \sum_{k = 1}^n \omega_k  {2n \choose n+k},
\end{equation} 
with $\omega_k$ repeating modulo 10 as given in Table \ref{tab.w10}, with $\omega_0 =2A_0$, 
$\omega_1 =A_1 - 2A_0$, and 
$\omega_2 = A_1 - 3A_0$. 

\begin{table}[H]
\centering
\begin{tabular}{c|ccccccccccc} 
  $k$ (mod 10) & 0&1 & 2 & 3& 4& 5& 6&7& 8 & 9 & 10\\ 
  \hline
  $\omega_k$ & $\omega_0$ & $\omega_1$ & $\omega_2$ & $-\omega_2$ & $-\omega_1$ & $-\omega_0$ & $-\omega_1$ & $-\omega_2$ & $\omega_2$ & $\omega_1$ & $\omega_0$
\end{tabular}
 \caption{Values of $\omega_k$ for signature $\{5,-5\}$}
\label{tab.w10}
\end{table}
\end{theorem}

\begin{proof}
Since the signature $\{5, -5\}$ appears in Table \ref{tab.1},    we  apply Theorem \ref{c.1} with  
$q=10$. The values for $\omega_0$ and $\omega_1$
are given to us in Table \ref{e.omegak012}. As for 
$\omega_2$, we learn from Table \ref{e.omegak012}
that $\omega_2 = A_2 -4A_1+2A_0$, but since $A_2 = 5A_1 - 5A_0$ then this becomes 
$\omega_2 = A_1 - 3A_0$. The other values follow from the conclusion of Theorem  
\ref{c.1} which tells us that 
    $\omega_k = \omega_{10 \pm k}$ and that $\omega_k = -\omega_{5-k}$. 
\end{proof}

We can now prove two nice identities for the binomial transforms of $F_{2n}$ and of $L_{2n}$, respectively. Both of these next two  corollaries feature an equality between a weighted sum of binomial coefficients (on the left) and a periodic weighted sum of binomial coefficients (on the right).

\begin{corollary}\label{c.12now}
    For $F_n$ the Fibonacci numbers, 
    then
\begin{align}
	 \sum_{i=0}^n {n \choose i} F_{2i} \label{e.12a}
&= \sum_{k=0}^n (-1)^{k+1}\left(\frac{k}{5}\right)  {2n \choose n+k} \\
	   \intertext{and also}
  5^n F_{2n} \label{e.12b}
&= \sum_{k = 0}^{2n} (-1)^{k+1}\left(\frac{k}{5}\right)  {4n \choose 2n+k},
	\end{align}
where   $\displaystyle {k \Leg 5}$ represents the Legendre symbol.
\end{corollary}

We note that both of the formulas in Corollary \ref{c.12now} are new, and equation (\ref{e.12b}) gives us the second part of equation (\ref{e.first}) from the opening paragraph of this paper. 

\begin{proof}
The binomial transform of $F_{2n}$ on the left of equation (\ref{e.12a}) is 
$0, 1, 5, 20, 75, \dots$, as
seen in \seqnum{A093131} where we also learn that it has signature $\{5, -5\}$. 
If we let $A_n$ equal that binomial transform, then we can apply 
Theorem \ref{t.5neg5}, and since $A_0 =0$ and $A_1 = 1$   we will get  $\omega_0 = 0$ and $\omega_1 = \omega_2 = 1$.
Thus, Table \ref{tab.w10} from
Theorem \ref{t.5neg5} becomes
\begin{center}
\begin{tabular}{c|ccccccccccc} 
  $k$ (mod 10) & 0&1 & 2 & \ \ 3& \ \ 4& 5& \ \ 6& \ \ 7& 8 & 9 & 10\\ 
  \hline
 $\omega_k$ &  0 & 1 & 1 & $-1$ & $-1$ & 0 & $-1$ & $-1$  & 1 & 1 & 0
\end{tabular}
\end{center}
and from this we see that $\omega_k = 
(-1)^{k+1} {k \Leg 5}$, giving us
equation (\ref{e.12a}).

As for  equation (\ref{e.12b}), we replace $n$ with $2n$ in equation (\ref{e.12a}) and 
we note from \seqnum{A093131} that 
$\sum_{i=0}^{2n} {2n \choose i} F_{2i}$ is equal to $5^n F_{2n}$.
\end{proof}

\begin{corollary}\label{c.12}
    For $L_n$ the Lucas numbers, 
    then
\begin{align}
   \sum_{i = 0}^n  {n \choose i} L_{2i}     \label{e.c13a}
&= 5\sum_{j \geq 0} (-1)^j {2n \choose n+5j}  
-  5 {2n-1 \choose n}  \qquad \mbox{for $n \geq 1$},\\
\intertext{and also}
    5^{n} L_{2n}    \label{e.c13b}
&= 5\sum_{j \geq 0} (-1)^{j} {4n \choose 2n+5j} 
-  5{4n-1 \choose 2n} \qquad \mbox{for $n \geq 1$}. 
\end{align}
\end{corollary}

We note that these two identities are not new; equation (\ref{e.c13a}) follows from 
Merca \cite[Corollary 8]{Merca}, and 
equation (\ref{e.c13b}) is just a special case of (\ref{e.c13a}). However, our approach is rather different from the
one in Merca's article.

\begin{proof}
If we let $A_n$ equal the binomial transform of $L_{2n}$ then the sequence $(A_n)_{n \geq 0}$ begins with 
$2, 5, 15, 50, 175, \dots$ and is
given by \seqnum{A020876} where we learn that it has signature $\{5, -5\}$. Hence, we can once again apply 
Theorem \ref{t.5neg5}, and since $A_0 =2$ and $A_1 = 5$ in this case,  we get  $\omega_0 = 4$, 
$\omega_1 = 1$, and $\omega_2 = -1$. 
Thus, Table \ref{tab.w10} from
Theorem \ref{t.5neg5} gives us the following table, where we have written $\omega_k$ as a sum of two periodic sequences in order to reveal a hidden pattern. 
\begin{center}
\begin{tabular}{c|cccccccccc} 
  $k$ (mod 10)  &1 &\ \  2 &  3& \ \ 4& \ \ 5& \ \ 6&  7& \ \ 8 & 9 & \ \ 10\\ 
  \hline
 $\omega_k$  & 1 & $-1$ & $1$ & $-1$ & $-4$ & $-1$ & $1$  & $-1$ & $1$ &\ \  4\\
 \hline
  first sequence &   1 & $-1$ & $1$ & $-1$ & \ \ $1$ & $-1$ & $1$  & $-1$ & $1$ & $-1$\\
  second sequence  & & & & & $-5$ & & & & & \ \ 5
\end{tabular}
\end{center}
As a result, Theorem \ref{t.5neg5} gives us
\begin{equation}\label{e.t11a}
\sum_{i=0}^n {n \choose i} L_{2i} = 2{2n \choose n}
-  
   \sum_{k=1}^n (-1)^{k} {2n \choose n+k} 
+  
   5\sum_{j\geq 1} (-1)^{j} {2n \choose n+5j}.
\end{equation}
If we adjust the right-hand side so that the two sums
each start at $0$ instead of $1$, then this becomes
\begin{equation}\label{e.midsum}
\sum_{i=0}^n {n \choose i} L_{2i} = -2{2n \choose n}
-  
   \sum_{k=0}^n (-1)^{k} {2n \choose n+k} 
+  
   5\sum_{j\geq 0} (-1)^{j} {2n \choose n+5j},
\end{equation}
and we do this because that middle sum is easily converted to a well-known identity. To be precise, if we replace ${2n \choose n+k}$ with ${2n \choose n-k}$, and then replace $k$ with $n-k$, then the middle sum in equation (\ref{e.midsum}) satisfies
\[
   \sum_{k=0}^n (-1)^{k} {2n \choose n + k} 
=  
  \sum_{k=0}^n (-1)^{k} {2n \choose n - k} 
= 
   \sum_{k=0}^n (-1)^{n-k} {2n \choose k},
\]
and this can be simplified further to 
\[
   (-1)^{n} \sum_{k=0}^n (-1)^{k} {2n \choose k}.
\]
From Benjamin and Quinn \cite[Identity 168]{BQ} we learn that the above sum is equal to 
${2n-1 \choose n}$, and so we can now write equation (\ref{e.t11a}) as
\begin{equation}\label{e.t11b}
\sum_{i=0}^n {n \choose i} L_{2i} = 
  - 2{2n \choose n}
-  
   {2n-1 \choose n} 
+  
   5\sum_{j \geq 0} (-1)^{j} {2n \choose n+5j}.
\end{equation}
It is relatively easy to show that 
$ - 2{2n \choose n} - {2n-1 \choose n} $ is 
equal to $-5 {2n-1 \choose n}$ so long as $n \geq 1$, and thus we obtain our desired equation (\ref{e.c13a}). 

As for  equation (\ref{e.c13b}), we replace $n$ with $2n$ in equation (\ref{e.c13a}) and 
we note from \seqnum{A020875} that 
$\sum_{i=0}^{2n} {2n \choose i} L_{2i}$ is equal to $5^n L_{2n}$. From this we can obtain equation (\ref{e.c13b}), as desired.
\end{proof}

As a result of Corollary \ref{c.12}, we can produce the following nice result about sums of every tenth binomial coefficient. As far as we can tell, this is a new formula. 

\begin{corollary}
    For $L_n$ the Lucas numbers and for $n \geq 1$, then 
\begin{equation}\label{e.28}
\sum_{j=-n}^n {4n \choose 2n+10j} = (2^{4n-1} +  L_{4n} +  5^nL_{2n})/5.
\end{equation}
\end{corollary}

\begin{proof}
    We begin with a result from Shibukawa \cite[equation (1.15)]{Shibukawa} with $r=2$, giving us
\[
L_{4n} = -2^{4n-1} + \frac{5}{2} \sum_{j=-n}^n {4n \choose 2n+5j}.
\]
We now re-write our equation (\ref{e.c13b}) from Corollary \ref{c.12} by using the identity
$-5{4n-1\choose 2n} = -\frac{5}{2}{4n \choose 2n}$; this gives us 
\[
 5^{n} L_{2n}   
= 
-\frac{5}{2}{4n \choose 2n} + 5\sum_{j = 0}^n (-1)^{j} {4n \choose 2n+5j} \qquad \mbox{for $n \geq 1$}, 
\]
and since each 
${4n \choose 2n+5j}$ is equal to 
${4n \choose 2n-5j}$, then we can re-write this again as 
\[
 5^{n} L_{2n}   
= 
\frac{5}{2}\sum_{j = -n}^n (-1)^{j} {4n \choose 2n+5j} \qquad \mbox{for $n \geq 1$}. 
\]
When we add this to Shibukawa's equation we get, after re-indexing our sum,
\[
L_{4n} + 5^n L_{2n} = -2^{4n-1} + 5\sum_{j = -n}^n  {4n \choose 2n+10j},
\]
and this gives us our desired equation (\ref{e.28}).
\end{proof}

\section{Proof of Theorem \ref{c.1}}\label{s.proof}

We now give the rather technical proof of our main result.

\begin{proof}[Proof of Theorem \ref{c.1}]
For $(A_n)_{n \geq 0}$ a linear
recurrence as given,  we let $f(x)$ represent its characteristic polynomial and we let $M$  be its order (equivalently, $M$ is the  degree of $f(x)$).
We are given  that this $f(x)$ is the minimal polynomial for $(2 \cos \pi/q)^2$ for some  integer $q\geq 3$, and since
$(2 \cos \pi/q)^2 = 2 \cos 2 \pi/q + 2$, then
$f(2x+2)$ is the minimal polynomial for $x = \cos 2 \pi/q$. It is well known \cite{Gurtas, Lehmer}
that $f(2x+2)$ has degree $M=\phi(q)/2$.

Furthermore, borrowing notation from  G\"{u}rta\c{s}  \cite{Gurtas}, if we set 
\[
S(q) = \{ p_i\, | \, \gcd(p_i,q) = 1 \mbox{\ \  and \ \ } 1 \leq p_i < q/2\}
\]
then the complete and distinct set of roots 
\cite{WZ} for $f(2x+2)$ are $\{\cos 2 \pi p_i/q \, | \, p_i \in S(q)\}$ and so likewise for $f(x)$ the roots are 
$\{(2\cos \pi p_i/q)^2 \, | \, p_i \in S(q)\}$. 
Now that we have the $M$ distinct roots of 
the characteristic polynomial $f(x)$, we know from Dubeau et.\ al.\ \cite[Theorem 1]{Dubeau} or Kelly and Peterson \cite[Theorem 3.7]{Kelley} that the Binet formula for $A_n$ is 
\begin{equation}\label{e.GenBinet}
A_n = \alpha_1 (2\cos \pi p_1/q)^{2n}  + \cdots + 
\alpha_M (2\cos \pi p_M/q)^{2n} = 
\sum_{i=1}^M \alpha_i (2 \cos \pi p_i/q)^{2n},
\end{equation}
with $p_1,  \dots p_M$ the distinct elements in the
set $S(q)$ and with $\alpha_1, \dots, \alpha_M$ determined by the initial values 
$A_0, A_1, A_2, \dots A_{M-1}$. 

We now call upon the well-known trigonometric power 
formula \cite{Weisstein} for cosine, 
\begin{equation}\label{e.cosine}
\cos^{2n} \theta = \frac{1}{2^{2n}} {2n \choose n} + \frac{1}{2^{2n-1}} \sum_{k=0}^{n-1} {2n \choose k} \cos 2(n-k) \theta.
\end{equation}
We simplify the above sum  by 
replacing $k$ with $n-k$ and noting that 
${2n \choose n-k} = {2n \choose n+k}$. After doing so, and then multiplying through by $2^{2n}$, we obtain
\begin{equation}\label{e.cosine2}
(2\cos \theta )^{2n} = {2n \choose n} +  \sum_{k=1}^{n} {2n \choose n+k} 
 2\cos 2 k \theta.
\end{equation}

Next, we apply equation (\ref{e.cosine2}) to each cosine in equation
(\ref{e.GenBinet}) to obtain 
\[
A_n =  \sum_{i=1}^M \alpha_i \left({2n \choose n} +  \sum_{k=1}^{n} {2n \choose n+k} 
 2\cos 2 k \pi p_i/q\right).
\]
We pull out the first term, and switch the order of summation on the rest, to get
\begin{equation}\label{e.Anproof}
A_n =  {2n \choose n} \sum_{i=1}^M \alpha_i  
     +  \sum_{k=1}^{n} \omega_k {2n \choose n+k}, 
\end{equation} 
where $\omega_k$ is  defined as
\begin{equation}\label{e.wkdef}
\omega_k = \sum_{i=1}^M \alpha_i \, 2 \cos 2k\pi p_i/q.
\end{equation}

 For the first term in equation (\ref{e.Anproof}), we note 
that if we take $n=0$ in equation 
(\ref{e.GenBinet}), we obtain
\[
A_0  = \sum_{i=1}^M \alpha_i,
\]
and so equation (\ref{e.Anproof}) becomes
\begin{equation}\label{e.Anproof2}
A_n =  A_0{2n \choose n} 
     +  \sum_{k=1}^{n} \omega_k {2n \choose n+k}, 
\end{equation} 
as desired. 

We now turn our attention to establishing equation (\ref{e.wk2}) for $\omega_k$. For this, we will use the 
Chebyshev polynomials $T_k(x)$ which have the useful property that $T_k(\cos \theta) = \cos k \theta$. We can write them as 
\begin{equation}\label{e.Tkdef}
T_k(x) = \sum_{j=0}^{\lfloor k/2 \rfloor} t_{k,j}x^{k-2j}
\end{equation}
such that our Chebyshev coefficients $t_{k,j}$, as seen in  the CRC handbook \cite[\S\,6.10.6]{CRC}, satisfy
\begin{equation}\label{e.tkj}
t_{k,j} = (-1)^j{k-j \choose j} 2^{k-2j-1}\frac{k}{k-j} 
\end{equation}
so long as $k-j \not= 0$. From our summation in equation (\ref{e.Tkdef}), we see that $0 \leq j \leq \lfloor k/2 \rfloor$, so the only case when $k-j=0$ is when both $j$ and $k$ are zero, and for that we define $t_{0,0} = 1$.

Next, we define $S_k(x)$ as
\begin{equation}
S_k(x) = 2 T_{2k}(x/2),
\end{equation}
which from equation (\ref{e.Tkdef}) gives us that 
\begin{equation}\label{e.Skdef}
S_k(x) = 2\sum_{j=0}^{k} t_{2k,j}
        \frac{x^{2k-2j}}{2^{2k-2j}}
\end{equation}
and furthermore we have that $S_k(2 \cos \theta) = 2 T_{2k}(\cos \theta)$ and thanks to the property of the Chebyshev polynomials, this equals $2\cos 2k\theta$. This means 
\[
S_k(2 \cos \pi p_i/q) = 2 \cos 2k \pi p_i/q,
\]
and so our formula for $\omega_k$ in equation (\ref{e.wkdef}) becomes
\[
\omega_k = \sum_{i=1}^M \alpha_i \cdot S_k(2 \cos \pi p_i/q).
\]
Thanks to equation (\ref{e.Skdef}) this becomes
\[ \omega_k  
    = \sum_{i=1}^M \alpha_i \cdot 2\sum_{j=0}^{k} \frac{t_{2k,j}}{2^{2k-2j}} (2 \cos \pi p_i/q)^{2k-2j}.
\]
We switch the order of summation to get
\begin{align*}
\omega_k &= 2 \sum_{j=0}^k \frac{t_{2k,j}}{2^{2k-2j}} \sum_{i=1}^M \alpha_i (2 \cos \pi p_i/q)^{2k-2j},\\
\intertext{and thanks to our formula for $A_n$ in equation (\ref{e.GenBinet}),
this becomes}
\omega_k &= 2 \sum_{j=0}^k \frac{t_{2k,j}}{2^{2k-2j}} A_{k-j}.
\end{align*}
We recall from our discussion above that $t_{0,0} = 1$ so this tells us that $\omega_0 = 2A_0$. For $k \geq 1$, we can use equation (\ref{e.tkj}) to replace the coefficients $t_{2k,j}$, and so the above equation becomes
\[
\omega_k =  \sum_{j=0}^k (-1)^j {2k-j \choose j} 
\frac{2k}{2k-j} A_{k-j} \qquad \mbox{for $k \geq 1$},
\]
as desired.

Next, since 
\[
\cos 2 (q \pm k) \pi p_i/q = \cos (2\pi p_i \pm 2 k\pi p_i/q) 
= \cos 2 k \pi p_i/q,
\] 
then from our definition of $\omega_k$ in equation (\ref{e.wkdef})
we have that $\omega_{k'} = \omega_{k}$ for $k' \equiv \pm k$ (mod $q$), as desired.

Finally, if $q$ is even, then for 
$k' \equiv q/2 - k$  {\rm (}mod $q${\rm )} we have
\[
\cos 2k'\pi p_i/q 
= 
\cos 2 (q/2 - k) \pi p_i/q 
= 
\cos ( \pi p_i + 2k\pi p_i/q) 
= 
-\cos 2 k \pi p_i/q,
\] 
where the last equality holds because since $q$ is even, then $p_i \in  S(q)$ must be odd. 
\end{proof}

\section{Conclusion}

We are rather surprised by how many different identities with weighted sums of binomial coefficients we could obtain from our Theorem \ref{c.1}, all thanks to the cosine formula in equation 
(\ref{e.cosine2}) for $(2 \cos \theta)^{2n}$. There is  a similar  formula for 
$(2 \cos \theta)^{2n+1}$, and we encourage the reader to explore the additional identities that this could provide. 

We might want to think about reversing the problem. Given a periodic sequence $(\omega_k)_{k \geq 0}$, what can we say about the periodic weighted sum
\[
A_n = \sum_{k=1}^n \omega_k {2n \choose n+k},
\]
and under what conditions on $\omega_k$ will 
$(A_n)_{n \geq 1}$ be a linear recurrence sequence? 

Also, while  Table \ref{tab.1}   only covers signatures of length $2$
or $3$, our Theorem \ref{c.1} applies to all signatures 
corresponding to  minimal polynomials for  $(2 \cos \pi/q)^2$ for all $q \geq 3$.  For  $q=11$, the signature 
$\{9,-28,35,-15,1\}$ has eight entries in the OEIS and it might be worthwhile to explore them. Another interesting example appears at $q=60$. 
The coefficients of the
minimal polynomial for 
$(2 \cos \pi/60)^2$ match up perfectly with the 
signature 
for the sequence \seqnum{A126569}
(which is related to the Cartan matrix for the Lie group $E_8$)
and so we could establish
a new identity for those terms as well.

Finally, we mention in closing another result on periodic binomial sums and the Fibonacci numbers. Moser \cite{Moser} gives this identity (adapted slightly to match our notation) for $F_{2n}$, which is a nice counterpoint to our equation (\ref{e.first}) at the beginning of the paper:
\[
F_{2n} = {2n-1 \choose n } - 2{2n-1 \choose n+2 } + {2n-1 \choose n+4 } 
- 2{2n-1 \choose n+6 } + \cdots.
\]

\section{Acknowledgements}
Our thanks to the editor and the referee for many helpful comments, and to Pioneer Academics for helping us to establish this research collaboration.

\bigskip
\hrule
\bigskip

\noindent 2010 {\it Mathematics Subject Classification}: Primary 11B39; 
Secondary 05A15, 05B45. 

\noindent \emph{Keywords: }  Fibonacci number, cosine, Catalan number, Pell number, binomial coefficient.

\bigskip
\hrule
\bigskip

\noindent (Concerned with sequences
\seqnum{A000108},
\seqnum{A000129},
\seqnum{A001075},
\seqnum{A001353},
\seqnum{A002203},
\seqnum{A007070},
\seqnum{A020875},
\seqnum{A020876},
\seqnum{A056236},
\seqnum{A080937},
\seqnum{A093131},
\seqnum{A126569},
\seqnum{A127677},
and
\seqnum{A359311}.)

\bigskip
\hrule
\bigskip

\vspace*{+.1in}
\noindent
Received October 10 2022;
revised version received March 11 2023. 
Published in {\it Journal of Integer Sequences}, Xxx xx 2023.

\bigskip
\hrule
\bigskip

\vskip .1in

\end{document}